\documentclass[10pt]{article}
\usepackage{ppfp}

\newcommand{\V}{\mathop{\rm V}}
\newcommand{\B}{\mathop{\rm B}}
\newcommand{\E}{\mathop{\rm E}}
\newcommand{\mv}{\mathop{\rm v.p.}}
\newcommand{\MISE}{\mathop{\rm MISE}}

\def\square{{\ \vrule height0.5em width0.5em depth-0.0em}}

{

\begin{document}

%% #1: title
%% #2: authors
%% #3: nr
%% #4: year
%% #5: name
%% #6: username
\thispagestyle{empty}
\begin{center}
    \Large\textrm{NORGES TEKNISK-NATURVITENSKAPELIGE}\\
    \textrm{UNIVERSITET}
\end{center}
\vspace*{2\baselineskip}
\begin{center}
    {\Large\textbf{Density Estimation Using the Sinc Kernel \\ }}
\end{center}
\begin{center}
    {\large by}
\end{center}
\begin{center}
    {\large\textrm{Ingrid K. Glad, Nils Lid Hjort and Nikolai Ushakov}}
\end{center}
\vspace*{2\baselineskip}
\begin{center}
    \Large
        \textrm{PREPRINT} \\
        \textrm{STATISTICS NO.~2/2007}
\end{center}
\vspace*{2\baselineskip}
\input epsf
\begin{center}
    \leavevmode
\end{center}
\vspace*{2\baselineskip}
\begin{center}
    \Large
    \textrm{NORWEGIAN UNIVERSITY OF SCIENCE AND}\\
    \textrm{TECHNOLOGY}\\
    \textrm{TRONDHEIM, NORWAY}
\end{center}
\vspace*{2\baselineskip}
\begin{center}
    \small
    This report has URL \texttt{http://www.math.ntnu.no/preprint/statistics/2007/S2-2007.ps}\\
    Nikolai Ushakov has homepage: \texttt{http://www.math.ntnu.no/$\sim$ushakov}\\
    E-mail: \texttt{ushakov@stat.ntnu.no}\\
    Address: Department of Mathematical Sciences, Norwegian University
    of Science and Technology, N-7491 Trondheim, Norway.
\end{center}
%%

%% #1: title
%% #2: authors
%% #3: number
%% #4: year
%% #5: your name
%% #6: your username

\begin{center}
{\huge \bf Density Estimation Using the Sinc Kernel}\\[0.5cm]

{\large
\begin{tabular}{c}
\large Ingrid K. Glad, Nils Lid Hjort\\ 
Department of Mathematics, University of Oslo, Norway\\
and\\
Nikolai Ushakov\\
Department of Mathematical Sciences \\
Norwegian University of Science and Technology \\
Trondheim, Norway
\end{tabular}
}
\end{center}

\begin{abstract}
This paper deals with the kernel density estimator based on the so-called sinc
(or Fourier integral) kernel $K(x)=(\pi x)^{-1}\sin x$. We study in detail both
asymptotic and finite sample properties of this estimator. It is shown that,
contrary to widespread opinion, the sinc estimator is superior to other
estimators in many respects: it is more accurate for quite moderate values of
the sample size, has better asymptotics in non-smooth case (the density to be 
estimated has only first derivative), is more convenient for the bandwidth
selection etc.
\end{abstract}

\noindent
{\em Key words}: Kernel estimation,
Sinc kernel, Fourier integral kernel, Mean integrated squared
error, Superkernels, Empirical characteristic function, Finite samples, Inequalities
\newpage

\noindent {\bf Contents}
\vskip 0.5cm

\noindent 1. Introduction

\noindent 2. The MISE of the estimator

\noindent 3. Comparison of the exact MISE of the sinc estimator and conventional
estimators (examples)

\noindent 4. Asymptotic superiority of the sinc estimator to conventional 
kernel estimators

\noindent 5. Comparison of the sinc estimator with superkernel estimators 

\noindent 6. Bandwidth selection

\noindent 7. Uniform consistency and estimation of the mode

\noindent 8. Inequalities

\noindent 9. Sinc estimator of derivatives

\vskip 0.8cm

\noindent {\bf 1. Introduction}
\vskip 0.5cm

Let $X_1,...,X_n$ be independent and
identically distributed random variables with the same probability density function
$f(x)$. We consider the problem of estimation of $f(x)$ nonparametrically. One of the
most popular methods is the kernel estimator
$$f_n(x;h)={1\over n}\sum_{i=1}^nK_h(x-X_i)\eqno(1.1)$$
where $K_h(x)=h^{-1}K(h^{-1}x)$, $K(x)$ is a 
kernel function (usually symmetric) and $h$
is the smoothing parameter (bandwith).

Typically, $K(x)$ is taken to be a probability density with at least
a couple of finite moments; this ensures that $f_n(x;h)$ itself
becomes a density function, and methods of Taylor expansion
and so on make it possible to analyse its behaviour to a
satisfactory degree. Recent monographs dealing with
general aspects of kernel density estimators include
Wand and Jones (1995) and Fan and Gijbels (1996).

However, the present paper deals with a non-standard choice
for $K(x)$, namely the so-called sinc kernel
$$K_s(x)={{\sin x}\over{\pi x}}$$
with the Fourier transform (characteristic function)
$$\psi_s(t)=\cases{
1&for $|t|\le 1$,\cr
0&for $|t|>1$.\cr}$$
The sinc kernel is not in $L_1$, that is its absolute value has 
infinite integral, but 
it is square integrable, and in
addition it is integrable in the sense
of the Cauchy principal value with 
$$\mv\int_{-\infty}^\infty K_s(x)dx=1,$$
in which
$$\mv\int_{-\infty}^\infty=\lim_{T\to\infty}
\lim_{\epsilon\to0}\left[\int_{-T}^{-\epsilon}+\int_\epsilon^T\right].$$
(In the following we will omit integration limits when the integral is to be taken over
the full real line.) Respectively, we have
$$\psi_s(t)=\mv\int e^{itx}K_s(x)dx.$$

Sometimes the sinc kernel is defined as $K(x)=\sin(\pi x)/(\pi x)$
with the Fourier transform 
$$\psi(t)=\cases{
1&for $|t|\le \pi$,\cr
0&for $|x|>\pi$.\cr}$$
Both functions $\sin x/(\pi x)$  and $\sin(\pi x)/(\pi x)$
integrate to one in the sense of the principal value, and 
the difference is only in the scale parameter.

The sinc kernel is a ``non-conventional" kernel, it takes negative values and is not 
integrable in the ordinary sense (we will say that a kernel is conventional if
it is a probability density function i.e. it is non-negative and integrates to one;
kernel estimators, based on conventional kernels, will be called conventional 
estimators). Respectively, realizations of the kernel estimator,
based on the sinc kernel (we will call it the sinc estimator; in some works, it is
called FIE --- Fourier integral estimator, see for example Davis, 1975 and Davis, 1977) 
are not probability 
density functions. This defect however can be easily corrected without loss of
the performance (see Glad et al. (2003)).

The sinc estimator has excellent asymptotic properties compared to conventional
kernel estimators when the density to be estimated is smooth i.e. has several
derivatives, see Davis (1975) and Davis (1977). If, for example, the density 
to be estimated is an analytic function of a certain
type, then the mean integrated squared error 
of the sinc estimator decreases as $n^{-1}$ as $n\to\infty$ while no one 
conventional estimator can provide the rate of convergence better than
$n^{-4/5}$. It is believed however that the performance of the sinc
estimator is good, roughly speaking, only for very large $n$, 
and only for very smooth $f(x)$. In addition, it is beleaved that
even for very large $n$ and very smooth $f(x)$, the sinc estimator is inferior 
to kernel estimators based on so-called superkernels --- non-conventional
kernels whose Fourier transform is continuous and equals 1 in some
neighbourhood of the origin. In this work, we try to show that these
beliefs are unjust.

In Section 3, we present examples 
demonstrating that the sinc estimator is more accurate than conventional
estimators for quite moderate values of the sample size. In Section 4,
we prove that the sinc estimator is asymptotically supereor to
conventional estimators (has a strictly better order of consistency)
even when $f(x)$ has only one derivative.
In Section 5, we make comparison of the sinc estimator with a 
superkernel estimator and show that the sinc estimator has better
properties, in particular, it is more accurate. 
In Section 6, we consider the problem of bandwidth
selection. This problem is solved easier and more effective for the sinc estimator
than for other estimators.
The problem of estimation of the mode is studied in
Section 7. Some useful inequalities for the MISE of the sinc
estimator are obtained in Section 8. The inequalities also show that
the performance of the sinc estimator is good not only for large
sample sizes but for moderate and even small too. In section 9, we study 
the problem of estimating derivatives. Here the sinc estimator is 
especially effective compared with other kernel estimators.

\vskip 1 cm

\noindent {\bf 2. The $\MISE$ of the estimator}
\vskip 0.5cm

Let $\hat f_n(x)$ be an estimator of $f(x)$ associated with the sample 
$X_1,...,X_n$. 
The customary performance criterion for density estimators
is the mean integrated squared error ($\MISE$), which is defined as
$$\MISE(\hat f_n(x))=\E \int[\hat f_n(x)-f(x)]^2dx.$$
In practice one seeks methods to minimise the $\MISE$ function.
The MISE is the sum of the integrated squared bias (denote it by
$\B(\hat f_n(x))$) and the integrated variance (denote it by
$\V(\hat f_n(x))$) of the estimator.

Denote the characteristic function of random variables $X_j$ by $\varphi(t)$
and the empirical characteristic function associated with the sample
$X_1,...,X_n$ by $\varphi_n(t)$:
$$\varphi(t)=\E e^{itX_j},\ \ 
\varphi_n(t)={1\over n}\sum_{j=1}^ne^{itX_j}.$$
Then the sinc estimator 
(in the rest of the paper, we denote it by $f_n(x;h)$) is
$$f_n(x;h)={1\over{\pi n}}\sum_{j=1}^n{{\sin[(x-X_j)/h]}\over{x-X_j}},$$
and its characteristic function equals 
$\varphi_n(t)\psi_s(ht)=\varphi_n(t)I_{[-1/h,1/h]}(t)$,
where, as usually, $I_A(t)$ denotes the indicator of the set $A$.

For a real valued function $g(x)$ we will use the following notation, 
provided the integrals exist:
$$\mu_k(g)=\int|x|^kg(x)dx,\ \ k=0,1,2,...,\ \ 
R(g)=\int g^2(x)dx.$$

The following lemma will be frequently used in the work.

{\bf Lemma 2.1.} \it For the sinc estimator,
$$\B(f_n(x;h))={1\over{2\pi}}\int_{|t|>1/h}|\varphi(t)|^2dt,\eqno(2.1)$$
$$\V(f_n(x;h))=
{1\over n}\cdot{1\over{2\pi}}\int_{-1/h}^{1/h}(1-|\varphi(t)|^2)dt,\eqno(2.2)$$
and
$$\MISE(f_n(x;h))={1\over{2\pi}}\int_{|t|>1/h}|\varphi(t)|^2dt+
{1\over n}\cdot{1\over{2\pi}}\int_{-1/h}^{1/h}(1-|\varphi(t)|^2)dt.\eqno(2.3)$$
where $\varphi(t)$ is the characteristic function of 
the density to be estimated. 
\rm

{\bf Corollary.} \it
$$\MISE(f_n(x;h))={1\over{\pi n h}}+R(f)-\left(1+{1\over n}\right){1\over\pi}
\int_0^{1/h}|\varphi(t)|^2dt.\eqno(2.4)$$
\rm

Equalities (2.1)--(2.4) can be found for example in Davis (1977).

Using Lemma 2.1, one can make some general remarks concerning the sinc estimator.
It is well known that in case of the kernel estimator with a 
conventional kernel, the necessary and sufficient condition of consistency
is $h\to0$ as $n\to\infty$ and $nh\to\infty$. In case of the sinc estimator, this
condition is also sufficient, but sometimes, not necessary. If the characteristic
function of the density to be estimated vanishes outside some interval containing
the origin: $\varphi(t)=0$ for $|t|>T$, then the necessary condition is milder:
the sinc estimator is consistent if $\limsup_{n\to\infty}h<1/T$. This circumstance
was pointed out in a number of works, see for example Davis (1977) or
Ibragimov and Khas'minskii (1982).

The second remark concerns the problem of selection of the smoothing parameter.
For a given $n$, minimum of the $\MISE$ function, as a function of $h$, can be
non-unique.
This means that there may exist several different optimal values of the
bandwidth $h$. This problem is considered more in details in Section 5.

\vskip 1 cm

\noindent {\bf 3. Comparison of the exact $\MISE$ of the sinc estimator and 
conventional estimators (examples)}
\vskip 0.5 cm

In this section, the exact $\MISE$ of the sinc estimator is compared with that 
of estimators based on some conventional kernels. 

{\it 3.1. Normal distribution.} Consider the standard normal density 
$$f(x)={1\over{\sqrt{2\pi}}}e^{-x^2/2}.$$
Let $f_n(x;h)$ be (as above) the sinc estimator of $f(x)$, and denote the kernel
estimator of $f(x)$, based on the normal kernel,
by $f_n^{(norm)}(x;h)$. In this subsection, we compare the performance of
$f_n(x;h)$ and $f_n^{(norm)}(x;h)$ for several finite values of the sample size.

The $\MISE$ of estimators $f_n(x;h)$ and $f_n^{(norm)}(x;h)$ is found explicitly:
$$\MISE(f_n)=
{1\over{\sqrt\pi}}
\left[
1-\left(1+{1\over n}\right)\Phi\left({{\sqrt2}\over{h}}\right)+
{1\over n}\left({1\over{h\sqrt\pi}}+{1\over2}\right)
\right]$$
and
$$\MISE(f_n^{(norm)})={1\over{2\sqrt\pi}}
\left(
1-2\sqrt{2\over{2+h^2}}+{1\over\sqrt{1+h^2}}+{1\over{nh}}-{1\over{n\sqrt{1+h^2}}}
\right),$$
where $\Phi(x)$ is the standard normal distribution function.
Values of $\inf_{h>0}\MISE(f_n)$ and $\inf_{h>0}\MISE(f_n^{(norm)})$
(and their ratio) are given in Table 1 
for the sample size $n=40$, $45$, $50$, $100$ and $1000$

%\vskip 1 cm

\centerline{Table 1}

\begin{center}
\begin{tabular}{|c|c|c|c|}
\hline
$n$&sinc&normal&sinc/normal\\
\hline
40&0.010141&0.01009&1.005\\
45&0.009203&0.009327&0.987\\
50&0.008436&0.00869&0.971\\
100&0.004699&0.005411&0.868\\
1000&0.000611&0.00103&0.593\\
\hline
\end{tabular}
\end{center}

%\vskip 1 cm

The Table shows that, under appropriate choice of the smoothing parameter for both
estimators, the sinc estimator is less accurate (but very little, only 0.5\%) than
the estimator, based on the normal kernel, when $n=40$, but it becomes better
already for $n=45$. For $n=100$, the sinc estimator is
about 15\% more accurate than normal. For large sample sizes ($>1000$), the sinc
estimator becomes several times better than normal (almost two times for $n=1000$).

One more advantage of the sinc estimator is that
$\MISE(f_n(x;h))\le\inf_{h>0}\MISE(f_n^{(norm)}(x;h))$
for a wide interval of values of the smoothing parameter $h$ ($0.4<h<0.53$ for $n=100$
and $0.25<h<0.46$ for $n=1000$). This means that, even if for the sinc estimator, $h$
is chosen quite far from its optimal value, the estimator is still better than the
normal estimator under the optimal choice of $h$.

{\it 3.2. Cauchy distribution.} Now consider the density 
$$f(x)={1\over{\pi(1+x^2)}}$$
(Cauchy distribution) with the characteristic function
$$\varphi(t)=e^{-|t|}.$$
Then the $\MISE$ of the sinc estimator is
$$\MISE(f_n)={1\over\pi}
\left[
{1\over2}\left(1+{1\over n}\right)e^{-2/h}+
{1\over n}\left({1\over h}-{1\over 2}\right)
\right].$$
For the sake of simplicity we consider the conventional estimator with the Cauchy
kernel
$$K(x)={1\over{\pi(1+x^2)}}.$$
Denote this estimator by
$f_n^{(Cauchy)}(x;h)$. Its $\MISE$ is found explicitly
$$\MISE(f_n^{(Cauchy)})={1\over\pi}
\left[
{1\over2}-{2\over{2+h}}+{1\over{2(1+h)}}+{1\over{2nh}}-{1\over{2n(1+h)}}
\right].$$
Results are presented in Table 2. They are very similar to those in the 
previous subsection.
Like in the normal case, the sinc estimator is superior to the conventional estimator
for $n\ge 45$ and becomes several times more accurate for $n>1000$.

%\vskip 1 cm

\centerline{Table 2}

\begin{center}
\begin{tabular}{|c|c|c|c|}
\hline
$n$&sinc&Cauchy&sinc/Cauchy\\
\hline
40&0.014776&0.014553&1.015\\
45&0.013542&0.01363&0.993\\
50&0.012516&0.012845&0.974\\
100&0.007346&0.00863&0.851\\
1000&0.0011&0.002126&0.517\\
\hline
\end{tabular}
\end{center}

\vskip 1 cm

\noindent {\bf 4. Asymptotic superiority of the sinc estimator to conventional
kernel estimators}
\vskip 0.5cm

In this section, we find conditions under which the sinc estimator 
$f_n(x;h)$ has a
strictly better order of consistency than conventional kernel estimators.
Let $K(x)$ be a square integrable conventional kernel (a square integrable 
probability density) and
$\hat f_n(x;h)$ --- the corresponding estimator.
It is known that if the density to be estimated $f(x)$ is 
at least three times differentiable, then 
$${{\MISE(f_n(x;h))}\over{\MISE(\hat f_n(x;h))}}\to0\ {\rm as}\ 
h\to0,\ nh\to\infty.\eqno(4.1)$$
It turns out that (4.1) holds under much broader conditions and
can take place even when $f(x)$ has only the first derivative.

Since $\MISE(f_n(x;h))=\MISE(\hat f_n(x;h))=\infty$ if $f(x)$ is
not square integrable, we, in the remainder of this section, suppose
that it is square integrable. In addition, without loss of generality,
we suppose that all conventional kernels under consideration are
symmetric.

Denote the characteristic function of $K(x)$ by $\psi(t)$.
For the integrated squared bias and integrated variance of $\hat f_n(x;h)$,
the following representations hold (see for example
Watson and Leadbetter, 1963)
$$\B(\hat f_n(x;h))={1\over{2\pi}}\int |\varphi(t)|^2(1-\psi(ht))^2dt,
\eqno(4.2)$$
$$\V(\hat f_n(x;h))={1\over n}\cdot{1\over{2\pi}}
\int (1-|\varphi(t)|^2)(\psi(ht))^2dt.
\eqno(4.3)$$

{\bf Lemma 4.1.} \it If
$$|\varphi(t)|=o\left({{1}\over{t^{5/2}}}\right)\ \ 
{\rm as}\ \ t\to\infty,\eqno(4.4)$$
then, under appropriate rescaling of the conventional kernel $K(x)$,
$$\B(f_n(x;h))=o(\B(\hat f_n(x;h)))\ \ {\rm as}\ \ h\to0.\eqno(4.5)$$
\rm

{\bf Proof.} There exist two positive numbers $\varepsilon$ and $c$
such that 
$$\psi(t)\le 1-\varepsilon t^2\ \ {\rm for}\ \ |t|\le c,\eqno(4.6)$$
see Loeve (1977). Without loss of generality one can suppose that
$c=1$ (otherwise a rescaling of $K(x)$ can be used). Rewrite (4.6) in
the form
$$(1-\psi(ht))^2\ge\varepsilon^2h^4t^4,\ \ |t|\le 1/h.\eqno(4.7)$$
Put $\lambda=1/h$,
$$R_1(\lambda)={1\over{\lambda^4}}\int_0^\lambda 
t^4|\varphi(t)|^2dt,$$
and
$$R_2(\lambda)={1\over{\lambda^4}}\int_\lambda^\infty 
|\varphi(t)|^2dt.$$
Note that (4.4) implies
$$|\varphi(\lambda)|^2=o\left({1\over{\lambda^5}}\int_0^\lambda 
t^4|\varphi(t)|^2dt
\right)\ \ {\rm as}\ \ \lambda\to\infty.\eqno(4.8)$$
Making use of (4.7), (4.8) and representations (2.1), (2.2), (4.2), and (4.3),
we obtain
$$\infty=\lim_{\lambda\to\infty}
\varepsilon^2\left[
{{{4\over{\lambda^5}}\int_0^\lambda 
t^4|\varphi(t)|^2dt}\over{|\varphi(\lambda)|^2}}
-1\right]=
\lim_{\lambda\to\infty}\varepsilon^2
{{R'_1(\lambda)}\over{R'_2(\lambda)}}=
\lim_{\lambda\to\infty}\varepsilon^2
{{R_1(\lambda)}\over{R_2(\lambda)}}=$$
$$=\lim_{h\to0}\varepsilon^2
{{h^4\int_0^{1/h}}t^4|\varphi(t)|^2dt
\over{\int_{1/h}^\infty|\varphi(t)|^2dt}}\le
\lim_{h\to0}
{{\int_0^{1/h}}|\varphi(t)|^2(1-\psi(ht))^2dt
\over{\int_{1/h}^\infty|\varphi(t)|^2dt}}\le$$
$$\le\lim_{h\to0}
{{\int_0^{\infty}}|\varphi(t)|^2(1-\psi(ht))^2dt
\over{\int_{1/h}^\infty|\varphi(t)|^2dt}}=
\lim_{h\to0}{{\B(\hat f_n(x;h))}\over{\B(f_n(x;h))}}$$
that implies (4.5).
\square

{\bf Theorem 4.1.} \it Let
$$|\varphi(t)|=o\left({{1}\over{t^{5/2}}}\right)\ \ 
{\rm as}\ \ t\to\infty.$$
Then
$$\inf_{h>0}\MISE(f_n(x;h))=
o(\inf_{h>0}\MISE(\hat f_n(x;h)))\ \ {\rm as}\ \ n\to\infty.\eqno(4.9)$$
\rm

{\bf Proof.} The main asymptotic term of the integrated variance of the
both estimators, conventional and sinc, has form $c/(nh)$, therefore,
due to Lemma 4.1 ($K(x)$ is scaled so that (4.5) holds)
$$\MISE(\hat f_n(x;h))\sim g(h)+{{c_1}\over{nh}}\eqno(4.10)$$
and
$$\MISE(f_n(x;h))\sim g(h)\varepsilon(h)+{{c_2}\over{nh}}\eqno(4.11)$$
where $g(h)\to0$ as $h\to0$, $\varepsilon(h)\to0$ as $h\to0$.
(4.10) and (4.11) evidently imply (4.9).
\square

Thus, under condition (4.4), the sinc estimator has strictly better
order of consistency than any conventional kernel. Condition (4.4)
can be satisfied even when $f(x)$ has only the first derivative.
For example, the density corresponding to the characteristic function
$(\sin t/t)^3$, that is the convolution of three standard uniform
densities, does not have the second detivative for $x=\pm 3$.
\vskip 0.8cm

\noindent {\bf 5. Comparison of the sinc estimator with superkernel
estimators}
\vskip 0.5cm

A superkernel is defined as a nonconventional kernel $K(x)$ whose
Fourier transform has form
$$\psi(t)=\cases{
1&for $|t|\le\Delta$,\cr
g(t)&for $\Delta\le|t|\le c\Delta$,\cr
0&for $|x|>c\Delta$,\cr}$$
where $g(t)$ is a real-valued , even function, satisfying the inequality
$|g(t)|\le1$ and chosen in such a way that $\psi(t)$ is continuous.
Superkernels were studied by Devroye (1992) in one-dimensional case and
by Politis and Romano (1999) in multidimensional case. The sinc kernel
can be considered as a limit case of superkernels as $c\to1$. Estimators, 
based on superkernels, have the same order of consistency as the sinc
estimator --- both are kernels of ``infinite order". Superkernels have one
advantage compared to the sinc kernel: corresponding estimates are integrable,
although this advantage is rather technical then essential. Advantages of the
sinc kernel compared to superkernels are simplicity (both for theoretical
analysis and practical use) and better solution of the problem of the
bandwidth selection. 

Politis and Romano (1999) state that ``$c=1$ is a
bad choice'', that is the sinc estimator is worse than superkernel
estimators.
In this section we compare accuracy of the sinc estimator 
with that of the most recommended by Politis and Romano (1999) superkernel 
estimator, namely, that for which $c=2$ and $g(t)$ is linear
on the interval $\Delta<t<2\Delta$. Thus
$$\psi(t)=\cases{
1&for $|t|\le\Delta$,\cr
2-{{|t|}\over{\Delta}}&for $\Delta\le|t|\le 2\Delta$,\cr
0&for $|t|>2\Delta$.\cr}$$
The estimator based on this superkernel is denoted in this section
by $\hat f_n(x;h)$. The sinc estimator is denoted, as earlier, by
$f_n(x;h)$.  The $\Delta$ does not play any role, therefore, for
convenience, we suppose that $\Delta<1$.

Now we make a comparison of the MISE of the sinc estimator and the
superkernel estimator under consideration for a broad class of
underlying distributions. We consider densities whose 
characteristic function $\varphi(t)$ satisfies condition
$$|\varphi(t)|^2={{a}\over{|t|^m}}\ \ {\rm for}\ \ |t|>c$$
where $a$, $c$ and $m$ are real positive constants, $m>3$.
Without loss of generality suppose that $a=1$. Then, for
sufficiently small $h$, namely, for $h<\Delta/c$, the integrated
squared bias of the both estimators is calculated explicitly.
Some elementary algebra shows that
$$\B(f_n(x;h))={{h^{m-1}}\over{\pi(m-1)}}\eqno(5.1)$$
and
$$\B(\hat f_n(x;h))={{h^{m-1}}\over{\pi\Delta^{m-1}}}
\left[
{{1}\over{m-1}}-
{{2(2^{m-2}-1)}\over{(m-2)2^{m-2}}}+
{{2^{m-3}-1}\over{(m-3)2^{m-3}}}
\right].\eqno(5.2)$$

For the integrated variance of the estimators the following 
equalities hold when $h<\Delta/c$ (also after some elementary algebra)
$$\pi n\V(f_n(x;h))={1\over h}-
\int_0^{\Delta/h}|\varphi(t)|^2dt-{{h^{m-1}}\over{\Delta^{m-1}}}
\left({{1-\Delta^{m-1}}\over{m-1}}\right)\eqno(5.3)$$
and
$$\pi n\V(\hat f_n(x;h))={{4\Delta}\over{3h}}-
\int_0^{\Delta/h}|\varphi(t)|^2dt-$$
$$-{{h^{m-1}}\over{\Delta^{m-1}}}
\left[
{{4(2^{m-1}-1)}\over{(m-1)2^{m-1}}}-
{{4(2^{m-2}-1)}\over{(m-2)2^{m-2}}}+
{{2^{m-3}-1}\over{(m-3)2^{m-3}}}
\right]
\eqno(5.4)$$

Once again, $\Delta$ does not play any role, so let us take it
to be equal $3/4$ (then the main asymptotic term of the integrated
variance of the superkernel estimator coincides with
that of the sinc estimator, and the comparison
becomes easier). For each $m\ge4$,
$${{4(2^{m-1}-1)}\over{(m-1)2^{m-1}}}-
{{4(2^{m-2}-1)}\over{(m-2)2^{m-2}}}+
{{2^{m-3}-1}\over{(m-3)2^{m-3}}}<{{1-(3/4)^{m-1}}\over{m-1}}
$$
and therefore the right hand side of (5.4) is greater than the 
right hand side of (5.3) that is the integrated variance of the 
sinc estimator is less than that of the superkernel estimator 
uniformly in $h$. 

Consider right hand sides of (5.1) and (5.2). For $m\le17$ 
$$\B(f_n(x;h))>\B(\hat f_n(x;h)),$$ 
while for $m\ge18$
$$\B(f_n(x;h))<\B(\hat f_n(x;h))$$
uniformly in $h$, and the ratio
$\B(f_n(x;h))/\B(\hat f_n(x;h))$ decreases with $m$
and tends to zero very fast as $m$ tends to infinity
(if, for example, $m>30$, then $\B(f_n(x;h))$ is more than ten
times smaller than $\B(\hat f_n(x;h))$). 

Thus, for more or less smooth underlying densities ($m=17$
approximately corresponds to the case when $f(x)$ has five
derivatives), the sinc estimator is more accurate than the
superkernel estimator under consideration. For really smooth
(many times differentiable) densities it is much more
accurate. Moreover, both the integrated variance and the
integrated squared bias of the sinc estimator are smaller
than those of the superkernel estimator {\it uniformly in
$h$}.
In non-smooth case (the fifth derivative of $f(x)$ does not 
exist), the accuracy of the two estimators is approximately
the same: the integrated variance of the sinc is smaller
for each $h$, while the integrated squared bias is greater
for each $h$. Now the estimators can be compared only under
condition that the bandwidth is chosen to be optimal for each 
of them.

Suming up and taking into account other advantages of the sinc
estimator, like simplicity, better solution of the bandwidth
selection etc., we must make conclusion that the sinc estimator
is preferable to the considered superkernel estimator.

\vskip 0.8cm

\noindent {\bf 6. Bandwidth selection}
\vskip 0.5cm

Representation of the $\MISE$, given by Lemma 1, suggests relatively simple 
rules of selection of the smoothing parameter $h$. Due to Corollary of Lemma 2.1,
$$\MISE(f_n(x;h))={1\over{\pi n h}}+R(f)-\left(1+{1\over n}\right){1\over\pi}
\int_0^{1/h}|\varphi(t)|^2dt.$$
Put $\delta=1/h$. Then
$${{\partial}\over{\partial\delta}}\MISE(f_n(x;h))={1\over{\pi n}}-
\left(1+{1\over n}\right){1\over\pi}|\varphi(\delta)|^2$$
and
$${{\partial^2}\over{\partial\delta^2}}\MISE(f_n(x;h))=-
\left(1+{1\over n}\right)
{1\over\pi}{{\partial}\over{\partial\delta}}|\varphi(\delta)|^2.$$
Therefore the optimal $\delta$ ($\delta$ minimizing the $\MISE$) must be a root
of the equation 
$$|\varphi(\delta)|={1\over{\sqrt{n+1}}},$$
and the optimal bandwidth $h_{\rm opt}$ is a solution of 
$$|\varphi(1/h)|={1\over{\sqrt{n+1}}}.\eqno(6.1)$$
\vskip 0.5cm

\noindent{\it 6.1 Normal rule}
\vskip 0.5cm

According to the normal scale rule, the bandwidth is selected so that it
minimizes the MISE if the underlying distribution is normal with the
variance $\sigma^2$, where unknown $\sigma^2$ is replaced by some its 
estimator. For the sinc estimator this rule is as follows. For a normal 
underlying distribution formula (6.1) becomes
$$e^{-{{\sigma^2}\over{2h^2}}}={{1}\over{\sqrt{n+1}}},$$
therefore the optimal value of $h$ (if $\sigma$ is known) is
$$h_{\rm opt}={{\sigma}\over{\ln(n+1)}},$$
and the normal rule bandwidth is
$$h_{\rm norm}={{\hat\sigma}\over{\ln(n+1)}}$$
where $\hat\sigma$ is some estimator of $\sigma$, for example,
the empirical standard deviation:
$$\hat\sigma=\left[{1\over n}\sum_{j=1}^n(X_j-\bar X)^2\right]^{1/2}$$
\vskip 0.5cm

\noindent{\it 6.2 Method based on the empirical characteristic function}
\vskip 0.5cm

Note that equation (6.1) always has solution 
(since $|\varphi(t)|\to0$ as $t\to\infty$)
but maybe non-unique. Consider all solutions of equation (6.1) for which
$|\varphi(\delta-0)|>1/\sqrt{n+1}$ and $|\varphi(\delta+0)|<1/\sqrt{n+1}$.
Denote them by $\delta_1,...,\delta_m$ 
(suppose for simplicity that $\delta_1<\delta_2<...<\delta_m$). 
Since $|\varphi(\delta)|$ decreases in some 
neighbourhood of each $\delta_i$,
$${{\partial}\over{\partial\delta}}|\varphi(\delta)|^2
\bigg|_{\delta=\delta_i}<0$$
and therefore
$${{\partial^2}\over{\partial\delta^2}}\MISE(f_n(x;h))
\bigg|_{\delta=\delta_i}>0.$$
Thus each $\delta_i$ is a local minimum of the $\MISE$.

The global minimum can be found by computing and comparison of the $\MISE$ at
$h=1/\delta_1,...,1/\delta_m$. This does not lead to large computational 
expenses because, if 
one uses (2.3) for the computation, the first integral in the right 
hand side of (2.3)
for $\delta=\delta_2,...,\delta_m$ is a part of this integral 
for $\delta=\delta_1$,
while the second integral for $\delta=\delta_1,...,\delta_{m-1}$ 
is a part of this 
integral for $\delta=\delta_m$.

The characteristic function $\varphi(t)$ however is (of course) 
unknown, therefore the
procedure, described above, is used to the empirical characteristic function
$\varphi_n(t)$ instead of $\varphi(t)$. Here one must take into account that
$\varphi_n(t)$ is an almost periodic function, and equation
$$|\varphi_n(\delta)|={1\over{\sqrt{n+1}}}$$
has infinitely many roots. Therefore, since
$$\lim_{n\to\infty}\sup_{|t|\le T_n}|\varphi_n(t)-\varphi(t)|=0$$
(see Cs\"org\H o and Totik, 1983),
where $T_n\to\infty$ and $\log T_n=o(n)$ as $n\to\infty$, it suffices to consider
only roots on the interval $[0,e^n]$. Of course it is not necessary
to calculate $\varphi_n(t)$ on such a wide interval: all roots of this interval
are contained in the interval $[0,\sqrt n]$, practically, 
in a much shorter interval.

\vskip 0.8cm

\noindent {\bf 7. Uniform consistency and estimation of the mode}
\vskip 0.5cm

In this section we prove that the sinc estimator
is uniformly consistent: it converges (in probability) to the true 
density function uniformly over the whole real line, and that the mode
of the sinc estimator is a consistent estimator of the mode of the true
density function. We formulate and prove results in the simplest form,
leaving possible generalizations to the reader.

Let $K(x)$ be a symmetric, differentiable, 
conventional kernel with finite variance
$\sigma^2$. Suppose also that its derivative has finite total variation
which we denote by $v$. Denote the characteristic function of $K(x)$
by $\psi(t)$ and the kernel estimator, based on $K(x)$, by 
$\hat f_n(x;h)$. As before, $f_n(x;h)$ denotes the sinc estimator.

{\bf Lemma 7.1.} \it Let the characteristic function $\varphi(t)$ of the
underlying distribution be integrable:
$$\int|\varphi(t)|dt<\infty.$$
Then
$$\sup_x|f_n(x;h)-\hat f_n(x;h)|\buildrel{a.s.}\over\longrightarrow 0\ 
{\rm as}\ n\to\infty, h\to0, nh\to\infty.$$
\rm

{\bf Proof.}
$$\sup_x|f_n(x;h)-\hat f_n(x;h)|\le$$
$$\le{{1}\over{2\pi}}
\bigg[
\int|\psi(ht)|\cdot|\varphi_n(t)-\varphi(t)|dt+
\int|\varphi(t)|\cdot|\psi(ht)-I_{[-1/h,1/h]}(t)|dt+$$
$$+\int_{-1/h}^{1/h}|\varphi_n(t)-\varphi(t)|dt
\bigg].$$
Now we prove that each of the three integrals in the square brackets tends
to zero as $n\to\infty, h\to0, nh\to\infty$. Denote these integrals by
$I_1$, $I_2$ and $I_3$, respectively. Then
$$I_1\le\int_{|t|\le n^2}|\varphi_n(t)-\varphi(t)|dt+
4\int_{n^2}^\infty|\psi(ht)|dt.$$
The first integral in the right hand side almost surely converges
to zero as $n\to\infty$ due to theorem 1 by Cs\"org\H o and Totik (1983).
To estimate the second integral, we use the inequality
$$|\psi(t)|\le{{v}\over{|t|^2}},$$
which holds for all $t$, see Ushakov and Ushakov (2000). Making use
of this inequality, we obtain
$$\int_{n^2}^\infty|\psi(ht)|dt\le{{v}\over{n^2h^2}}\to0
\ {\rm as}\ nh\to\infty.$$
So, $I_1\buildrel{a.s.}\over\longrightarrow 0$ as $n\to\infty, nh\to\infty$.
$$I_2\le2\int_0^{1/\sqrt h}[1-\psi(ht)]dt+
4\int_{1/\sqrt h}^\infty|\varphi(t)|dt.$$
The second integral tends to zero as $h\to0$ because $|\varphi(t)|$
is integrable. To estimate the first integral, we use the inequality
$$\psi(t)\ge 1-{{\sigma^2t^2}\over{2}}$$
which holds for all $t$, see Ushakov (1999).
Making use of this inequality, we obtain
$$\int_0^{1/\sqrt h}[1-\psi(ht)]dt\le{{\sigma^2}\over{6}}\sqrt h\to0
\ {\rm as}\ h\to0.$$
Thus $I_2\to0$ as $h\to0$.

Finally, if $nh\to\infty$, then $1/h\le cn$ with some constant $c$,
and therefore
$$I_3\le\int_{|t|\le cn}|\varphi_n(t)-\varphi(t)|dt
\buildrel{a.s.}\over\longrightarrow 0\ 
{\rm as}\ n\to\infty$$
due to the mensioned theorem by Cs\"org\H o and Totik (1983).
\square

{\bf Remark.} The condition of the lemma ($\varphi$ is integrable) implies
that $f(x)$ is uniformly continuous but is a little more restrictive.
It is satisfied for example when $f(x)$ is differentiable and $f'(x)$
has finite total variation.

{\bf Theorem 7.1.} \it Let $\varphi(t)$ be integrable.
Then
$$\sup_x|f_n(x;h)-f(x)|\buildrel{P}\over\longrightarrow 0\ 
{\rm as}\ n\to\infty, h\to0, nh^2\to\infty.\eqno(7.1)$$
\rm

{\bf Proof.} Let $K(x)$ be any conventional kernel satisfying 
conditions of both Lemma 7.1
of this section and Theorem 3A  of Parzen (1962). Then, due to Theorem 3A
by Parzen (1962),
$$\sup_x|\hat f_n(x;h)-f(x)|\buildrel{P}\over\longrightarrow 0\ 
{\rm as}\ n\to\infty, h\to0, nh^2\to\infty,\eqno(7.2)$$
and, due to
Lemma 7.1,
$$\sup_x|f_n(x;h)-\hat f_n(x;h)|\buildrel{P}\over\longrightarrow 0\ 
{\rm as}\ n\to\infty, h\to0, nh^2\to\infty.\eqno(7.3)$$
(7.2) and (7.3) evidently imply (7.1).
\square

Denote a mode of $f(x)$ by $\theta$. Suppose it is unique. Let $\theta_n$
be a mode of the sinc estimate.

{\bf Theorem 7.2.} \it Let $\varphi(t)$ be integrable.
Then
$$\theta_n\buildrel{P}\over\longrightarrow\theta\ 
{\rm as}\ n\to\infty, h\to0, nh^2\to\infty.$$
\rm

Proof of the theorem coincides with that of the second part of Theorem 3A
by Parzen (1962).

\vskip 0.8cm

\noindent {\bf 8. Inequalities}
\vskip 0.5cm

In this section, we derive some upper bounds for the $\MISE$ of the sinc estimator. 
In addition to practical importance (evaluation of the sample size sufficient for
achieving a given accurancy etc.), these inequalities throw more light onto
properties of the sinc estimator (especially for finite samples) 
and make easier comparison of this estimator
with other estimators. 

We define the 0-th derivative of a function as
the function itself: $f^{(0)}(x)=f(x)$ (as usually).
We need below the following form of the Parseval equality:
let $f(x)$ be $m$ times differentiable probability density
function ($m\ge0$), its $m$-th derivative $f^{(m)}(x)$ be square integrable,
and $\varphi(t)$ be the corresponding characteristic function. Then
$$\int(f^{(m)}(x))^2dx={1\over{2\pi}}\int t^{2m}|\varphi(t)|^2dt.\eqno(8.1)$$

{\bf Theorem 8.1.} \it Let $f(x)$ be $m$ times differentiable ($m\ge0$), and
its $m$-th derivative is square integrable. Then
$$\MISE(f_n(x;h))<\varepsilon(h)h^{2m}R(f^{(m)})+{1\over{\pi nh}}\eqno(8.2)$$
where $\varepsilon(h)\le1$ for all $h$ and $\varepsilon(h)\to0$ as $h\to0$.
\rm

{\bf Proof.} Estimate the first summand in the right hand side of (2.3). 
Making use of (8.1), we obtain
$${1\over{2\pi}}\int_{|t|>1/h}|\varphi(t)|^2dt=
h^{2m}{1\over{2\pi}}\int_{|t|>1/h}(1/h)^{2m}|\varphi(t)|^2dt\le$$
$$\le h^{2m}{1\over{2\pi}}\int_{|t|>1/h}t^{2m}|\varphi(t)|^2dt=
h^{2m}{1\over{2\pi}}\int t^{2m}|\varphi(t)|^2dt-$$
$$-h^{2m}{1\over{2\pi}}\int_{-1/h}^{1/h} t^{2m}|\varphi(t)|^2dt
=h^{2m}{1\over{2\pi}}\int t^{2m}|\varphi(t)|^2dt
\left(1-{{\int_{-1/h}^{1/h} 
t^{2m}|\varphi(t)|^2dt}\over{\int t^{2m}|\varphi(t)|^2dt}}\right)=$$
$$=\varepsilon(h)h^{2m}\int(f^{(m)}(x))^2dx=\varepsilon(h)h^{2m}R(f^{(m)}),$$
where
$$\varepsilon(h)=
1-{{\int_{-1/h}^{1/h} 
t^{2m}|\varphi(t)|^2dt}\over{\int t^{2m}|\varphi(t)|^2dt}}$$
evidently satisfies conditions of the theorem: 
$\varepsilon(h)\le1$ and $\varepsilon(h)\to0$ as $h\to0$.

For the second summand in the right hand side of (2.3) we have
$${1\over n}\cdot{1\over{2\pi}}\int_{-1/h}^{1/h}(1-|\varphi(t)|^2)dt<
{1\over n}\cdot{1\over{2\pi}}\int_{-1/h}^{1/h}dt={1\over{\pi nh}}.$$
Thus we finally obtain (8.2). \square

{\bf Corollary 1.} \it Let conditions of theorem 1 be satisfied. Then 
$$\MISE(f_n(x;h))<h^{2m}R(f^{(m)})+{1\over{\pi nh}}.\eqno(8.3)$$
\rm

Putting in (8.3)
$$h=\left[{{1}\over{2\pi nmR(f^{(m)})}}\right]^{1\over{2m+1}}$$
(this $h$ minimizes the right hand side of (8.3)), we get

{\bf Corollary 2.} \it Let conditions of Theorem 8.1 be satisfied. Then
$$\inf_{h>0}\MISE(f_n(x;h))<{{1+2m}\over{(2\pi m)^{(2m)/(2m+1)}}}
R(f^{(m)})^{1/(2m+1)}n^{-{{2m}\over{2m+1}}}.$$
\rm 

{\bf Corollary 3.} \it Let conditions of Theorem 8.1 be satisfied. Then
$$\inf_{h>0}\MISE(f_n(x;h))=o\left(n^{-2m/(2m+1)}\right),\ \ \ n\to\infty.$$
\rm

{\bf Corollary 4.} \it If $f(x)$ is two times differentiable, and its second
derivative is square integrable, then
$$\inf_{h>0}\MISE(f_n(x;h))=o(n^{-4/5}),\ \ n\to\infty$$
and
$$\inf_{h>0}\MISE(f_n(x;h))<{5\over{4\pi}}(4\pi R(f''))^{1/5}n^{-4/5}.$$
\rm

To obtain more sensitive and accurate estimates we use 
one more characteristic of the density
to be estimated --- its total variation
(or the total variation of its derivatives). For a function $g(x)$, we denote 
its total variation by $Vr(g)$.
In the general case, two conditions, 
finiteness of the total variation and square
integrability, are not comparable: 
a function may have finite variation but be not 
square integrable and vice versa. But for densities, 
square integrability is milder
then finiteness of the total variation: 
a density, having bounded variation, is square 
integrable. Thus finiteness of the total variation of a density 
(or its derivatives) 
is a little more restrictive condition than its square integrability. But
on the other hand, assumption that some derivative of the density to be
estimated has finite total variation (together with the use of the sinc kernel)
allows one to improve the order of decreasing  
the $\MISE$. For example, if the density
is two times differentiable and $Vr(f'')<\infty$, 
then $\MISE(f_n(x;h))=O(n^{-5/6})$,
$n\to\infty$, instead of $o(n^{-4/5})$ (see below).

{\bf Theorem 8.2.} \it Let $f(x)$ be $m$ times differentiable ($m\ge0$), 
and its $m-$th derivative has finite total variation.
Then
$$\MISE(f_n(x;h))\le 
h^{2m+1}{{Vr(f^{(m)})^2}\over{(2m+1)\pi}}+{{1}\over{\pi nh}}.\eqno(8.4)$$
\rm

{\bf Proof.} For all $t$,
$$|\varphi(t)|\le{{Vr(f^{(m)})}\over{|t|^{m+1}}},$$
see Ushakov and Ushakov (2000). Making use of 
this inequality, estimate the first summand in the
right hand side of (2.3):
$${1\over{2\pi}}\int_{|t|>1/h}|\varphi(t)|^2dt\le
{{Vr(f^{(m)})^2}\over{\pi}}\int_{1/h}^\infty{{dt}\over{t^{2m+2}}}
={{Vr(f^{(m)})^2}\over{(2m+1)\pi}}h^{2m+1}.$$
For the second summand of the right hand side of (2.3) 
we have (see proof of Theorem 8.1)
$${1\over n}\cdot{1\over{2\pi}}\int_{-1/h}^{1/h}(1-|\varphi(t)|^2)dt<
{1\over{\pi nh}}.$$
So, we obtain (8.4). \square

{\bf Corollary.} \it Let conditions of Theorem 8.2 be satisfied. Then
$$\inf_{h>0}\MISE(f_n(x;h))\le
{{2(m+1)}\over{(2m+1)\pi}}Vr(f^{(m)})^{1/(m+1)}n^{-(2m+1)/(2m+2)}.$$
\rm

For example, if $m=2$, we get
$$\inf_{h>0}\MISE(f_n(x;h))\le{6\over{5\pi}}Vr(f'')^{1/3}n^{-5/6}.$$

Following Watson and Leadbetter (1963) and
Davis (1975), we will say that a characteristic function $\varphi(t)$
decreases exponentially with degree $\alpha$ and 
coefficient $\rho$ ($\rho>0$, $0<\alpha\le2$)
if
$$|\varphi(t)|\le Ae^{-\rho|t|^\alpha}\eqno(8.5)$$
for some constant $A$. Davis (1975, Theorem 4.1) proved
that if the characteristic function of the density to be 
estimated satisfies (8.5),
then 
$$\lim_{n\to\infty}he^{\rho/h^\alpha}|\B(f_n(x;h))|=0$$
The next theorem makes this more precise.

{\bf Theorem 8.3.} \it Let 
$${1\over{2\pi}}\int e^{\rho|t|^\alpha}|\varphi(t)|^2dt
=C<\infty.$$
Then
$$\MISE(f_n(x;h))\le\varepsilon(h)Ce^{-\rho/h^\alpha}+
{1\over{\pi nh}},\eqno(8.6)$$
where $0<\varepsilon(h)<1$ and $\varepsilon(h)\to0$ as $h\to0$.
\rm

{\bf Proof} is similar to that of Theorem 8.1. 
For the first summand of the right hand
side of (2.3), we have
$${1\over{2\pi}}\int_{|t|>1/h}|\varphi(t)|^2dt<
e^{-\rho/h^\alpha}{1\over{2\pi}}
\int_{|t|>1/h}e^{\rho|t|^\alpha}|\varphi(t)|^2dt=$$
$$=\varepsilon(h)Ce^{-\rho/h^r}$$
where
$$\varepsilon(h)=
1-{{\int_{-1/h}^{1/h}e^{\rho|t|^\alpha}|\varphi(t)|^2dt}\over
{\int e^{\rho|t|^\alpha}|\varphi(t)|^2dt}}.$$ 
\square

It is difficult to find explicitly $h$ minimizing the right hand side of (8.6), 
therefore
we take $h$ for which the right hand side of (8.6) has a simple form. 
Namely, put
$$h=\left({1\over\rho}\ln n\right)^{-1/\alpha},$$
then
$$\MISE(f_n(x;h))<
\left(C+{{(\ln n)^{1/\alpha}}\over{\pi\rho^{1/\alpha}}}\right){1\over n}<
\left(C+{1\over{\pi\rho^{1/\alpha}}}\right){{(\ln n)^{1/\alpha}}\over{n}}$$
provided $n>2$. If, for example, $f(x)$ is the standard normal density, then
$$\MISE(f_n(x;h))<
\left({1\over{\sqrt{2\pi}}}+{{\sqrt 2}\over\pi}\right){{\sqrt{\ln n}}\over{n}}$$

{\bf Corollary.} \it 
Let the characteristic function $\varphi(t)$ of $f(x)$ decreases
exponentially with degree $r$ and coefficient $\rho$ (i.e. satisfies (8.5)). 
Then
$$\lim_{n\to\infty}e^{c\rho/h^\alpha}|\B(f_n(x;h))|=0$$
for any $c<2$.
\rm

\vskip 0.8 cm

\noindent {\bf 9. Estimation of derivatives}
\vskip 0.5cm

The sinc estimator is especially superior to 
conventional estimators when a derivative
of $f(x)$ is estimated. Let $f(x)$ be $r$ times differentiable. 
Suppose that one needs to 
estimate the $r-$th derivative $f^{(r)}(x)$. 
A natural way is to estimate $f^{(r)}(x)$
by the $r-$th derivative of a kernel estimator of $f(x)$, 
provided that the kernel is
$r$ times differentiable. So, let
$$\hat f_n(x;h)={1\over{nh}}\sum_{j=1}^nK\left({{x-X_j}\over{h}}\right)$$
be a kernel estimator of $f(x)$, and $K^{(r)}(x)$ exists. Then the estimator
$$\hat f_n^{(r)}(x;h)=
{1\over{nh^{r+1}}}\sum_{j=1}^nK^{(r)}\left({{x-X_j}\over{h}}\right)$$ 
is used for estimation of $f^{(r)}(x)$.

If $K(x)$ is a conventional kernel, then the $\MISE$ of the estimator 
$\hat f_n^{(r)}(x;h)$
is represented in the form (provided that $f(x)$ has $r+2$ derivatives 
and $K(x)$ has finite
second moment)
$$\MISE(\hat f_n^{(r)}(x;h))={1\over4}h^4\mu_2(K)^2R(f^{(r+2)})+
{1\over{nh^{2r+1}}}R(K^{(r)})+$$
$$+o\left(h^4+{1\over{nh^{2r+1}}}\right),\ \ h\to0.$$
Therefore the optimal $\MISE$ is of order $n^{-4/(2r+5)}$, 
i.e. the rate becomes slower 
for higher values of $r$, and the difficulty increases, 
see Wand and Jones (1995)
and Stone (1982). 
In this section, we show that the sinc estimator $f_n^{(r)}(x;h)$
is free of this difficulty and the estimator $f_n^{(r)}(x;h)$ of $f^{(r)}(x)$ has
almost the same order of consistency as
the estimator $f_n(x;h)$ of $f(x)$ 
even for large values of $r$ (of course, if the density to be estimated is smooth enough).
First we formulate an analog of Lemma 2.1 for derivatives.

{\bf Lemma 9.1.} \it For the sinc estimator,
$$\MISE(f_n^{(r)}(x;h))={1\over{2\pi}}\int_{|t|>1/h}t^{2r}|\varphi(t)|^2dt+
{1\over n}\cdot{1\over{2\pi}}\int_{-1/h}^{1/h}t^{2r}(1-|\varphi(t)|^2)dt.
\eqno(9.1)$$
where $\varphi(t)$ is the characteristic function of 
the density to be estimated. 
\rm

{\bf Proof.} Due to the Parseval-Plancherel identity we have
$$\MISE(f_n^{(r)}(x;h))=\E\int(f_n^{(r)}(x;h)-f^{(r)}(x))^2dx=$$
$$={1\over{2\pi}}\E\int t^{2r}|\varphi_n(t)I_{[-1/h,1/h]}(t)-\varphi(t)|^2dt=$$
$$={1\over{2\pi}}\int\E\Big[(\varphi_n(t)I_{[-1/h,1/h]}(t)-\varphi(t))
(\overline{\varphi_n(t)}\ I_{[-1/h,1/h]}(t)-\overline{\varphi(t)})\Big]dt=$$
$$={1\over{2\pi}}\int_{-1/h}^{1/h} t^{2r}\E|\varphi_n(t)|^2dt-$$
$$-{1\over{2\pi}}\int_{-1/h}^{1/h} t^{2r}
\Big(\varphi(t)\E\overline{\varphi_n(t)}+
\overline{\varphi(t)}\E\varphi_n(t)\Big)dt+
{1\over{2\pi}}\int t^{2r}|\varphi(t)|^2dt.\eqno(9.2)$$
It is easy to see that
$$\E\varphi_n(t)=\varphi(t),\eqno(9.3)$$
$$\E\overline{\varphi_n(t)}=\overline{\varphi(t)}\eqno(9.4)$$
and
$$\E|\varphi_n(t)|^2=\E\left|{1\over n}\sum_{j=1}^ne^{itX_j}\right|^2
=\E\left[{1\over n}\sum_{j=1}^ne^{itX_j}\cdot{1\over n}
\sum_{k=1}^ne^{-itX_k}\right]=$$
$${1\over{n^2}}\left[n+\sum_{j\not=k}e^{it(X_j-X_k)}\right]=
{1\over n}+\left(1-{1\over n}\right)|\varphi(t)|^2.
\eqno(9.5)$$
Substituting (9.3)---(9.5) to the right hand side of (9.2), we obtain (9.1). 
\square

Making use of Lemma 9.1 we obtain following
analogs of theorems of Section 8.

{\bf Theorem 9.1.} \it Let $f(x)$ be $r+m$ times differentiable ($r,m\ge0$), and
its $r+m$-th derivative is square integrable. Then
$$\MISE(f_n^{(r)}(x;h))<
\varepsilon(h)h^{2m}R(f^{(r+m)})+{1\over{\pi(2r+1)nh^{2r+1}}}$$
where $\varepsilon(h)\le1$ for all $h$ and $\varepsilon(h)\to0$ as $h\to0$.
\rm

{\bf Corollary.} \it Let conditions of Theorem 9.1 be satisfied. Then
$$\MISE_{h>0}(f_n^{(r)}(x;h))
\le C_{m,r}R(f^{(r+m)})^{(2r+1)/(2r+2m+1)}n^{-2m/(2r+2m+1)}$$
where
$$C_{m,r}=(2\pi m)^{-2m/(2r+2m+1)}+
{{(2\pi m)^{(2r+1)/(2r+2m+1)}}\over{\pi(2r+1)}}.$$
\rm

{\bf Theorem 9.2.} \it Let 
$$
{1\over{2\pi}}\int t^{2r}e^{\rho|t|^\alpha}|\varphi(t)|^2dt=C<\infty.$$
Then
$$\MISE(f_n^{(r)}(x;h))\le\varepsilon(h)
Ce^{-\rho/h^\alpha}+{1\over{\pi nh^{2r+1}}},$$
where $0<\varepsilon(h)<1$ and $\varepsilon(h)\to0$ as $h\to0$.
\rm

{\bf Corollary.} \it Let conditions of Theorem 9.2 be satisfied. Then
$$\inf_{h>0}\MISE(f_n^{(r)}(x;h))<
\left(C+{{(\ln n)^{(2r+1)/\alpha}}\over{\pi\rho^{(2r+1)/\alpha}}}\right)
{1\over n}.$$
\rm

{\bf Theorem 9.3.} \it Let the characteristic function $\varphi(t)$ of
$f(x)$ satisfies the condition: there exists $T>0$ such that
$f(t)=0$ for $|t|>T$. Then, if
$$h\le {1\over T},$$
then
$$\MISE(f_n^{(r)}(x;h))\le{1\over{\pi nh^{2r+1}}}.$$
In particular, if $h=const=1/T$, then
$$\MISE(f_n^{(r)}(x;h))\le{T^{2r+1} \over{\pi n}}.$$
\rm

Proofs of Theorems 9.1 -- 9.3 are similar to those of 
theorems of Section 8 therefore we leave them to the reader.
\vskip 1 cm

\noindent {\bf References}
\vskip 0.5cm

Cs\"org\H o, S. and Totik, V. (1983). On how long interval is the empirical
characteristic function uniformly consistent? {\it Acta Sci. Math.}, {\bf 45},
141-149.

Davis, K.B. (1975). Mean square error properties of density estimates. 
{\it Ann. Statist.}, {\bf 3}, no. 4, 1025-1030.

Davis, K.B. (1977). Mean integrated square error properties of density estimates. 
{\it Ann. Statist.}, {\bf 5}, no. 3, 530-535.

Devroye, L. (1992). A note on the usefulness of superkernels in density 
estimates. {\it Ann. Statist.}, {\bf 20}, 2037-2056.

Fan, J., Gijbels, I. (1996). {\it Local polynomial modelling and its applications.}
Monograps on Statistics and Applied Probability. Chapman and Hall, London.

Glad, I.K., Hjort, N.L. and Ushakov, N.G. (2003). 
Correction of density estimators
that are not densities. {\it Scand. J. Statist.}, {\bf 30}, no. 2, 415-427.

Ibragimov, I.A., Khas'minskii, R.Z. (1982). Estimation of 
distribution density belonging to a class of entire functions.
{\it Theory Probab. Applic.}, 
{\bf 27}, No. 3, 551-562.

Loeve, M. (1977). {\it Probability Theory.} Springer, Berlin.

Parzen, E. (1962). On estimation of a probability density function
and its mode. {\it Ann. Math. Statist.}, {\bf 33}, 1065-1076.

Politis, D.N., Romano, J.P. (1999). Multivariate density estimation with 
general flat-top Kernels of infinite order. {J. Multiv. Anal.},
{\bf 68}, 1-25.

Stone, C.J. (1982). 
Optimal global rates of convergence on nonparametric regression.
{\it Ann. Statist.}, {\bf 10}, no. 4, 1040-1053.

Ushakov, N.G. (1999). 
{\it Selected Topics in Characteristic Functiond}. VSP, Utrecht.

Ushakov, V.G. and Ushakov, N.G. (2000). Some inequalities for characteristic
functions of densities with bounded variation. 
{\it Moscow Univ. Comput. Math. Cybernet.,} 
no. 3, 45-52.

Wand, M.P. and Jones, M.C. (1995). {\it Kernel smoothing.} Chapman and Hall, 
London.

Watson, G.S. and Leadbetter, M.R. (1963) On the estimation of 
the probability density, I. {\it Ann. Math. Statist.,} Vol. 34, 
480-491.

\end{document}